\documentclass[a4paper,10pt]{article}
\usepackage{stmaryrd}
\usepackage{amsfonts}
\usepackage{bbm}
\usepackage{amscd}
\usepackage{mathrsfs}
\usepackage{latexsym,amssymb,amsmath,amscd,amscd,amsthm,amsxtra}
\usepackage[dvips]{graphicx}
\usepackage[utf8]{inputenc}
\usepackage[T1]{fontenc}
\usepackage{lmodern}
\usepackage{amssymb}
\usepackage[all]{xy}
\usepackage{nicefrac,mathtools,enumitem}
\usepackage{microtype}
\usepackage{xcolor}
\newtheorem{thm}{Theorem}[section]
\newtheorem{lem}[thm]{Lemma}

\newtheorem{pro}[thm]{Proposition}
\newtheorem{ex}[thm]{Example}
\newtheorem{rmk}[thm]{Remark}
\newtheorem{defi}[thm]{Definition}

\newcommand {\emptycomment}[1]{}

\newcommand{\Sym}{\mathrm {Sym}}
\newcommand{\lon }{\,\rightarrow\,}
\newcommand{\be }{\begin{equation}}
\newcommand{\ee }{\end{equation}}

\newcommand{\pf}{\noindent{\bf Proof.}\ }

\newcommand{\g}{\frkg}

\newcommand{\huaV}{\mathcal{V}}
\newcommand{\huaW}{\mathcal{W}}

\newcommand{\frkg}{\mathfrak g}

\def\qed{\hfill ~\vrule height6pt width6pt depth0pt}

\newcommand{\br}[1]{   [ \cdot,    \cdot  ]_\frkg   }

\newcommand{\Id}{\rm{Id}}

\newcommand{\dM}{\mathrm{d}}

\newcommand{\HF}{\mathsf{F}}

\newcommand{\Ad}{\mathsf{Ad}}

\newcommand{\gl}{\mathfrak {gl}}

\newcommand{\so}{\mathfrak {so}}

\newcommand{\ad}{\mathsf{ad}}

\newcommand{\Img}{\mathrm{Im}}

\begin{document}
\title{
{ Generalized metric $n$-Leibniz algebras and generalized orthogonal representation of metric Lie algebras}
\thanks
 {
Research supported by NSFC (11471139) and NSF of Jilin Province (20170101050JC).
 }
}
\author{Lina Song  and   Rong Tang\\
Department of Mathematics, Jilin University,\\\vspace{3mm}
 Changchun 130012, Jilin, China
\\
Email:songln@jlu.edu.cn, tangrong16@mails.jlu.edu.cn}

\date{}


\maketitle
\begin{abstract}
We introduce the notion of a generalized metric $n$-Leibniz algebra and show that there is a one-to-one correspondence between generalized metric $n$-Leibniz algebras and faithful generalized orthogonal  representations of metric Lie algebras (called Lie triple data). We further show that there is also a one-to-one correspondence between generalized orthogonal derivations (resp. generalized orthogonal automorphisms) on  generalized metric $n$-Leibniz algebras and   Lie triple datas.
\end{abstract}

\tableofcontents

\section{Introduction}
Ternary Lie algebras (3-Lie algebras) or more generally $n$-ary Lie algebras are  natural generalization of Lie  algebras.  They were introduced and studied by  Filippov in  \cite{Filippov}, and can be traced back to Nambu (\cite{Nambu}). See \cite{repKasymov,Ling,Liu-Makhlouf-Sheng,ST} and the review article \cite{review} for more details. This type of algebras appeared also in the algebraic formulation of Nambu Mechanics \cite{Nambu}, generalizing Hamiltonian mechanics by considering two hamiltonians, see \cite{Gautheron,T}. Moreover, 3-Lie algebras appeared in String Theory and
M-theory. In \cite{Basu},
 Basu and Harvey suggested to replace the Lie algebra appearing in the Nahm equation by a 3-Lie algebra for the  lifted Nahm equations.
Furthermore, in the context of Bagger-Lambert-Gustavsson model of multiple
M2-branes, Bagger-Lambert  managed to construct, using a ternary bracket,  an $N=2$
 supersymmetric version of the worldvolume theory of the M-theory membrane, see \cite{BL0}. These metric 3-Leibniz algebras (generalized 3-Lie algebras) have many applications, see \cite{Ch,CS,PJE,Escobar} for more details. Metric 3-Lie algebras and metric $n$-Lie algebras were further studied in \cite{BWL,N,ZhangBai}.

The notion of an $n$-Leibniz algebra was introduced in \cite{Loday-n-Leibniz} as a generalization of an $n$-Lie algebra and a Leibniz algebra \cite{Loday,Loday and Pirashvili}. See also \cite{deI} for more results.
 Through fundamental objects one may   represent    an $n$-Leibniz algebra by  a Leibniz algebra \cite{DT}.  Motivated by the work in \cite{DFMP}, where the authors established a one-to-one correspondence between metric 3-Leibniz algebras and faithful orthogonal representation of metric Lie algebras, it is natural to investigate the $n$-ary case. However, for the usual metric $n$-Leibniz algebras, where $n>3$, one can not use the method provided in \cite{DFMP}. We overcome this difficulty by introducing the notion of a generalized metric $n$-Leibniz algebra, where the ``metric'' is a symmetric non-degenerate $(n-1)$-linear form satisfying some compatibility conditions. We also introduce the notion of a generalized orthogonal representation of a Lie algebra and show that there is a one-to-one correspondence between generalized metric $n$-Leibniz algebras and faithful generalized orthogonal representation of metric Lie algebras. We also lift this one-to-one correspondence to the level of generalized orthogonal derivations and generalized orthogonal automorphisms.

The paper is organized as follows. In Section 2, we give a review of $n$-Leibniz algebras and metric Lie algebras. In Section 3, we construct a faithful generalized orthogonal representation of a metric Lie algebra from a generalized metric $n$-Leibniz algebra. In Section 4, we construct a generalized metric $n$-Leibniz algebra from a faithful generalized orthogonal representation of a metric Lie algebra. In Section 5, we show that there is a one-to-one   correspondence between generalized orthogonal derivations  on  generalized metric $n$-Leibniz algebras and   Lie triple datas. In Section 6, we show that there is a one-to-one   correspondence between   generalized orthogonal automorphisms  on  generalized metric $n$-Leibniz algebras and   Lie triple datas.

\vspace{2mm}

 \noindent {\bf Acknowledgement:} We give our warmest thanks to Yunhe Sheng  for very helpful discussions. We also thank the referees for very helpful comments that improve the paper.

\section{Preliminaries}
In this paper, we work over real field $\mathbb R$ and all the vector
spaces are finite-dimensional.

 \begin{defi}{\rm(\cite{Loday-n-Leibniz})}
   An $n$-Leibniz algebra is a vector space $\huaV$ equipped with an $n$-linear map $[\cdot,\cdots,\cdot]:\huaV\times\cdots\times\huaV\lon\huaV$ such that for all $u_1,\cdots,u_{n-1}$, $v_1,\cdots, v_n\in \huaV$, the following {\bf fundamental} identity  holds:
\begin{eqnarray}\label{fundamental-identity}
[u_1,\cdots,u_{n-1},[v_1,\cdots, v_n]]=\sum_{i=1}^n[v_1,\cdots,v_{i-1},
     [u_1,\cdots,u_{n-1},v_i],v_{i+1},\cdots,v_{n}].
\end{eqnarray}
 \end{defi}

In particular, if $n=2$, we obtain the notion of a Leibniz algebra \cite{Loday,Loday and Pirashvili}. If the $n$-linear map $[\cdot,\cdots,\cdot]$ is skew-symmetric, we obtain the notion of an $n$-Lie algebra \cite{Filippov}.

\begin{defi}{\rm(\cite{Loday-n-Leibniz})}
  A {\bf derivation} on an $n$-Leibniz algebra $(\huaV,[\cdot,\cdots,\cdot])$ is a linear map $\dM_{\huaV}\in\gl(\huaV)$, such that for all $u_1,\cdots,u_{n}\in \huaV$ the following equality holds:
\begin{eqnarray}
\dM_{\huaV}[u_1,\cdots,u_{n}]=\sum_{i=1}^n[u_1,\cdots,u_{i-1},\dM_{\huaV} u_i,u_{i+1},\cdots,u_n].
\end{eqnarray}
\end{defi}

Define $D:\otimes^{n-1}\huaV\lon\gl(\huaV)$
by
\begin{eqnarray}\label{defi:D}
D(u_1,\cdots,u_{n-1})u_n=[u_1,\cdots,u_{n-1},u_n],\,\,\,\,\forall u_1,\cdots,u_{n-1},u_n\in\huaV.
\end{eqnarray}
Then the fundamental identity \eqref{fundamental-identity} is the  condition that $D(u_1,\cdots,u_{n-1})$ is a derivation on the $n$-Leibniz algebra $(\huaV,[\cdot,\cdots,\cdot])$.

On $\otimes^{n-1}\huaV$, one can define a new bracket operation $[\cdot,\cdot]_{\HF}$ by
\begin{equation}\label{leibniz-bracket}
  [U,V]_{\HF}=\sum_{i=1}^{n-1}v_1\otimes\cdots\otimes v_{i-1}\otimes[u_1,\cdots,u_{n-1},v_i]\otimes v_{i+1}\otimes\cdots\otimes v_{n-1},
\end{equation}
for all $U=u_1\otimes\cdots \otimes u_{n-1},~V=v_1\otimes\cdots \otimes v_{n-1}\in \otimes^{n-1}\huaV$.
It is proved in \cite{DT} that $(\otimes^{n-1}\huaV,[\cdot,\cdot]_{\HF})$ is a  Leibniz algebra.  The fundamental identity \eqref{fundamental-identity} is equivalent to
\begin{eqnarray}\label{leibniz-morphism}
[D(U),D(V)]=D([U,V]_{\HF}).
\end{eqnarray}
Thus, we obtain that $D$ is a Leibniz algebra homomorphism from $\otimes^{n-1}\huaV$ to $\gl(\huaV)$.

\begin{defi}{\rm (\cite[Definition 2]{BH})}
Let $(A,\cdot)$ be a nonassociative algebra and $\omega$ a non-degenerate symmetric bilinear form on $A$.
\begin{itemize}
\item[\rm(i)] If $\omega(x\cdot y,z)=\omega(x,y\cdot z)$, then we say that $\omega$ is {\bf associative-invariant};
\item[\rm(ii)]  If $\omega(x\cdot y,z)=-\omega(y,x\cdot z)$, then we say that $\omega$ is {\bf(left) ad-invariant};
  \item[\rm(iii)]If $\omega(x\cdot y,z)=-\omega(x,z\cdot y)$, then we say that $\omega$ is {\bf(right) ad-invariant}.
\end{itemize}
\end{defi}

A non-degenerate symmetric bilinear form $\omega$ satisfies at least two of the
preceding definitions if and only if $(A,\cdot)$ is an anticommutative algebra. Since a Lie bracket  is skew-symmetric, we obtain that left ad-invariant, right ad-invariant and associative-invariant non-degenerate symmetric bilinear forms on a Lie algebra are the same. See \cite{BH} for more details.

Recall that
a   Lie algebra $(\g,[\cdot,\cdot])$ is a said to be {\bf metric} if it is equipped with
a symmetric non-degenerate   bilinear form $\omega$ which is {\bf(left) ad-invariant}, that is:
\begin{eqnarray}\label{ad-invariant}
\omega([x,y],z)=-\omega(y,[x,z]),\quad \forall x,y,z\in\g.
\end{eqnarray}

Moreover, there is a natural notion of orthogonal derivations and automorphisms on  metric Lie algebras.
\begin{defi}
Let $(\g,[\cdot,\cdot],\omega)$ be a metric Lie algebra. A derivation $\dM_{\g}$ on the Lie algebra $(\g,[\cdot,\cdot])$ is  called {\bf orthogonal} if the following equality holds:
\begin{eqnarray}\label{o-derivation}
\omega(\dM_{\g}x,y)+\omega(x,\dM_{\g}y)=0.
\end{eqnarray}
\end{defi}

\begin{defi}
Let $(\g,[\cdot,\cdot],\omega)$ be a metric Lie algebra. An automorphism $\Phi_{\g}$ on the Lie algebra $(\g,[\cdot,\cdot])$ is  called {\bf orthogonal} if the following equality holds:
\begin{eqnarray}\label{o-automorphism}
\omega(\Phi_{\g}x,\Phi_{\g}y)=\omega(x,y).
\end{eqnarray}
\end{defi}

\section{Construction of a Lie triple data from a generalized metric $n$-Leibniz algebra}

Let $\huaV$ be a vector space and $\huaV^*$ its dual space. Denote by $\Sym^k(\huaV^*)$ the vector space of symmetric tensors of order $k$ on $\huaV^*$. Any $\phi\in\Sym^k(\huaV^*)$ induces a linear map $\phi^\sharp:\huaV\longrightarrow \Sym^{k-1}(\huaV^*)$ by
$$
\phi^\sharp(u)(v_1,\cdots,v_{k-1})=\phi(u,v_1,\cdots,v_{k-1}),\quad\forall u,v_1,\cdots,v_{k-1}\in\huaV.
$$
$\phi\in\Sym^k(\huaV^*)$  is said to be non-degenerate if the induced map $\phi^\sharp:\huaV\longrightarrow \Sym^{k-1}(\huaV^*)$ is non-degenerate, that is, $\phi^\sharp(u)=0$ if and only if $u=0.$

\begin{defi}\label{n-leibniz-algebra}
A {\bf generalized metric $n$-Leibniz algebra} is an $n$-Leibniz algebra $(\huaV,[\cdot,\cdots,\cdot])$ equipped with a symmetric non-degenerate $(n-1)$-tensor $S\in\Sym^{n-1}(\huaV^*)$   satisfying the following axioms for all $u_1,\cdots,u_{n-1}$, $v_1,\cdots, v_{n-1}\in \huaV:$
\begin{itemize}
        %\item[\rm(a)] the {\bf nondegeneracy} of $S$ in the first entry, that is,
            %\begin{eqnarray}\label{nondegeneracy}
%S(u,v_1,\cdots,v_{n-2})=0\,\,\,\,\Rightarrow u=0;
            %\end{eqnarray}
        \item[\rm(a)] the {\bf unitarity} condition
            \begin{eqnarray}\label{unitarity}
          \sum_{i=1}^{n-1}S(v_1,\cdots, v_{i-1},[u_1,\cdots,u_{n-1},v_i],v_{i+1},\cdots, v_{n-1})=0;
            \end{eqnarray}
        \item[\rm(b)] the {\bf symmetry} condition
             \begin{eqnarray}\label{symmetry}
             S([u_1,u_2,\cdots,u_{n-1},v_1],v_2,\cdots, v_{n-1})=S([v_1,\cdots, v_{n-1},u_1],u_2,\cdots,u_{n-1}).
             \end{eqnarray}
        \end{itemize}
\end{defi}

We denote a generalized   metric $n$-Leibniz algebra by $(\huaV,[\cdot,\cdots,\cdot],S)$.

\emptycomment{
\begin{rmk}
When $n=2$, a generalized metric $2$-Leibniz algebra defined above is a Leibniz algebra $(\huaV,[\cdot,\cdot])$ equipped with a linear function $S\in\huaV^*$ such that the following condition holds:
\begin{eqnarray}\label{Flinear-forms}
S([u,v])=0,\quad \forall  u,v\in \huaV.
\end{eqnarray}
\end{rmk}
}

\begin{rmk}
When $n=3$ in Definition \ref{n-leibniz-algebra}, we obtain the notion of a generalized metric $3$-Leibniz algebra, which is the same as the generalized metric Lie $3$-algebra introduced in \cite[Definition 1]{DFMP}. See \cite{DFMP} for more applications of generalized metric Lie $3$-algebras in the BLG theory.
\end{rmk}

\begin{pro}
Let $(\huaV,[\cdot,\cdots,\cdot],S)$ be a  generalized   metric $n$-Leibniz algebra.
Then we have
\begin{eqnarray*}
\sum_{i=1}^{n-1}[v_i,v_1,\cdots,v_{i-1},\hat{v_i},v_{i+1},\cdots,v_{n-1},v_n]=0,\quad \forall v_1,\cdots,v_{n}\in\huaV.
\end{eqnarray*}
\end{pro}

\pf For all $v_1,\cdots,v_n,u_1,\cdots,u_{n-2}$, we have
\begin{eqnarray*}
&&S([v_1,v_2,\cdots,v_{n-1},v_n],u_1,\cdots,u_{n-2})\\&\stackrel{\eqref{symmetry}}{=}&S([v_n,u_1,\cdots,u_{n-2},v_1],v_2,\cdots,v_{n-1})\\
                                                  &\stackrel{\eqref{unitarity}}{=}&-\sum_{i=2}^{n-1}S([v_n,u_1,\cdots,u_{n-2},v_i],v_1,v_2,\cdots,v_{i-1},\hat{v_i},v_{i+1},\cdots,v_{n-1})\\
                                                  &\stackrel{\eqref{symmetry}}{=}&-\sum_{i=2}^{n-1}S([v_i,v_1,v_2,\cdots,v_{i-1},\hat{v_i},v_{i+1},\cdots,v_{n-1},v_n],u_1,\cdots,u_{n-2}).
\end{eqnarray*}
Since $S$ is non-degenerate, we have
\begin{eqnarray*}
[v_1,v_2,\cdots,v_{n-1},v_n]=-\sum_{i=2}^{n-1}[v_i,v_1,\cdots,v_{i-1},\hat{v_i},v_{i+1},\cdots,v_{n-1},v_n],
\end{eqnarray*}
which finishes the proof. \qed

\begin{defi}\label{orthogonal-derivation}
Let $(\huaV,[\cdot,\cdots,\cdot],S)$ be a generalized metric $n$-Leibniz algebra. A derivation $\dM_{\huaV}$ on the $n$-Leibniz algebra $(\huaV,[\cdot,\cdots,\cdot])$ is called {\bf generalized orthogonal} if the following equality holds:
\begin{eqnarray}\label{orthogonal-derivation-eq}
\sum_{i=1}^{n-1}S(v_1,\cdots,\dM_{\huaV} v_i,\cdots,v_{n-1})=0,
\end{eqnarray}
for all $v_1,\cdots,v_{n-1}\in\huaV.$
\end{defi}

\begin{defi}\label{orthogonal-automorphism}
Let $(\huaV,[\cdot,\cdots,\cdot],S)$ be a generalized metric $n$-Leibniz algebra. An automorphism $\Phi_{\huaV}$ on the $n$-Leibniz algebra $(\huaV,[\cdot,\cdots,\cdot])$ is called {\bf generalized orthogonal} if the following equality holds:
\begin{eqnarray}\label{orthogonal-automorphism-eq}
S(\Phi_{\huaV}v_1,\cdots,\Phi_{\huaV}v_{n-1})=S(v_1,\cdots,v_{n-1}),
\end{eqnarray}
for all $v_1,\cdots,v_{n-1}\in\huaV.$
\end{defi}

\begin{defi}\label{orthogonal-representation}
Let $\g$ be a Lie algebra and $\huaV$ a vector space equipped with a symmetric non-degenerate $(n-1)$-tensor $S\in\Sym^{n-1}(\huaV^*)$. A representation $\rho:\g\lon\gl(\huaV)$  is called {\bf generalized orthogonal} if the following equality holds:
\begin{eqnarray}\label{orthogonal-representation-invariant}
\sum_{i=1}^{n-1}S(w_1,\cdots,w_{i-1},\rho(x)w_i,w_{i+1},\cdots w_{n-1})=0,
\end{eqnarray}
for all $x\in\g$ and $w_1,w_2,\cdots,w_{n-1}\in \huaV$.
\end{defi}

We denote a generalized orthogonal representation by $(\rho,\huaV,S)$. When $n=3$, we recover the usual notion of an orthogonal representation of a Lie algebra.

We introduce the notion of a Lie triple data, which is the main object in this paper.

\begin{defi}
  A {\bf Lie triple data} consists of the following structure:
  \begin{itemize}
    \item[\rm(i)] a metric Lie algebra $(\g,[\cdot,\cdot],\omega)$;
    \item[\rm(ii)] a vector space $\huaV$ equipped with a symmetric non-degenerate $(n-1)$-tensor $S\in\Sym^{n-1}(\huaV^*)$;
    \item[\rm(iii)] a faithful generalized orthogonal representation  $\rho:\g\longrightarrow\gl(\huaV)$.
  \end{itemize}
  We will denote a Lie triple data by $(\g,\huaV,\rho)$.
\end{defi}

\emptycomment{
\begin{ex}\label{example-3}
The $sl_2(\mathbb C)$ is the
three-dimensional  complex simple Lie algebra generated by
elements $e,f$ and $h$ with Lie brackets
\begin{eqnarray*}
[e,h]=2e,\quad [f,h]=-2f,\quad  [e,f]=h.
\end{eqnarray*}
For all $n\ge 0$, there exists one irreducible representation $(V_n;\rho_n)$ of $sl_2(\mathbb C)$ (up to isomorphism) of dimension $n+1$.  The vector space $V_n$ with a basis $\{v_0,v_1,\cdots,v_n\}$ and the representation of $sl_2(\mathbb C)$ on $V_n$ is given by
\begin{eqnarray*}
\rho_n(h)v_i=(n-2i)v_i,\quad\rho_n(e)v_i=(n-i+1)v_{i-1},\quad\rho_n(f)v_i=(i+1)v_{i+1},
\end{eqnarray*}
here we set $v_{-1}=v_{n+1}=0.$ Work out the average-operators with respect to representation $(V_n;\rho_n)$ will be sophisticated unless we get help from the computers. However, for $n=1$
$$
   Tv_0=r_{11}e+r_{21}f+r_{31}h,\quad Tv_1=r_{12}e+r_{22}f+r_{32}h.
   $$
It is straightforward to see that $T$ is an average-operator with respect to representation $(V_1;\rho_1)$ if and only if
$$
?
$$
For $n=2$, we have
$$
   Tv_0=r_{11}e+r_{21}f+r_{31}h,\quad Tv_1=r_{12}e+r_{22}f+r_{32}h,\quad Tv_2=r_{13}e+r_{23}f+r_{33}h.
$$
It is straightforward to see that $T$ is an average-operator with respect to representation $(V_2;\rho_2)$\footnote{Since $sl_2(\mathbb C)$ is a complex simple Lie algebra. The representation $(V_2;\rho_2)$ is isomorphic to the adjoin representation $(sl_2(\mathbb C);\ad)$.} if and only if
$$
?
$$
\end{ex}
}

\subsection{From an $n$-algebra to a Lie algebra}
Let $(\huaV,[\cdot,\cdots,\cdot],S)$ be a generalized metric $n$-Leibniz algebra. Let $\g=\Img D\subset\gl(\huaV)$, where $D$ is given by \eqref{defi:D}.

\begin{pro}
 $(\g,[\cdot,\cdot]_C)$ is a Lie subalgebra of $\gl(\huaV)$, where $[\cdot,\cdot]_C$ denote the commutator Lie bracket on $\gl(\huaV)$.
\end{pro}

\pf By the fundamental identity \eqref{fundamental-identity}, we have
\begin{eqnarray*}
&&D(u_1,\cdots,u_{n-1})(D(v_1,\cdots,v_{n-1})v_n)-D(v_1,\cdots,v_{n-1})(D(u_1,\cdots,u_{n-1})v_n)\\
&=&\sum_{i=1}^{n-1}D(v_1,\cdots,v_{i-1},D(u_1,\cdots,u_{n-1})v_i,v_{i+1},\cdots,v_{n-1})v_n.
\end{eqnarray*}
Hence, we have
\begin{eqnarray}\label{commutator}
[D(u_1,\cdots,u_{n-1}),D(v_1,\cdots,v_{n-1})]_C=\sum_{i=1}^{n-1}D(v_1,\cdots,v_{i-1},D(u_1,\cdots,u_{n-1})v_i,v_{i+1},\cdots,v_{n-1})\in\g,
\end{eqnarray}
which shows that $[\g,\g]_C\subset\g$. The proof is finished. \qed\vspace{3mm}

Furthermore, we claim that $\g$ is a metric Lie algebra, that is, there is a symmetric non-degenerate   $\ad$-invariant bilinear form $\omega$ on $\g$. Actually, this bilinear form is defined by\footnote{By $D(u_1,\cdots,u_{n-1})=0$, for all $v\in \huaV$, we have $$D(u_1,\cdots,u_{n-1})v=[u_1,\cdots,u_{n-1},v]=0.$$ Thus, the definition of $\omega$ is
well-defined.}
\begin{eqnarray}\label{bilinear-form}
\omega(D(u_1,\cdots,u_{n-1}),D(v_1,\cdots,v_{n-1}))=S(D(u_1,\cdots,u_{n-1})v_1,v_2,\cdots,v_{n-1}).
\end{eqnarray}

\begin{pro}\label{pro:metLie}
The bilinear form $\omega$ on $\g$ defined by \eqref{bilinear-form} is symmetric, non-degenerate and $\ad$-invariant. Consequently, $(\g,\omega)$ is a metric Lie algebra.
\end{pro}

\pf By the symmetry condition \eqref{symmetry} of a generalized metric $n$-Leibniz algebra, we have
\begin{eqnarray*}
\omega(D(u_1,\cdots,u_{n-1}),D(v_1,\cdots,v_{n-1}))&=&S(D(u_1,\cdots,u_{n-1})v_1,v_2,\cdots,v_{n-1})\\
                                             &=&S([u_1,\cdots,u_{n-1},v_1],v_2,\cdots,v_{n-1})\\
                                             &=&S([v_1,,v_2,\cdots,v_{n-1},u_1],u_2,\cdots,u_{n-1})\\
                                             &=&\omega(D(v_1,\cdots,v_{n-1}),D(u_1,\cdots,u_{n-1})).
\end{eqnarray*}
Thus, the bilinear form $\omega$ is symmetric.

To prove non-degeneracy, let $x\in\g\subset\gl(\huaV)$ be such that $\omega(x,D(u_1,\cdots,u_{n-1}))=0$ for all $u_1,\cdots,u_{n-1}\in\huaV$. Thus, we have
$$S(x(u_1),u_2,\cdots,u_{n-1})=0.$$
By the nondegeneracy of $S$, we have $x(u_1)=0$ for all $u_1\in\huaV$, which implies that $x=0$.

Finally, we prove the ad-invariance of the bilinear form $\omega$:
\begin{eqnarray*}
&&\omega(D(u_1,\cdots,u_{n-1}),[D(v_1,\cdots,v_{n-1}),D(w_1,\cdots,w_{n-1})]_C)\\
&\stackrel{\eqref{commutator}}{=}&\omega(D(u_1,\cdots,u_{n-1}),\sum_{i=1}^{n-1}D(w_1,\cdots,w_{i-1},D(v_1,\cdots,v_{n-1})w_i,w_{i+1},\cdots,w_{n-1}))\\
&=&\sum_{i=1}^{n-1}\omega(D(u_1,\cdots,u_{n-1}),D(w_1,\cdots,w_{i-1},D(v_1,\cdots,v_{n-1})w_i,w_{i+1},\cdots,w_{n-1}))\\
&\stackrel{\eqref{bilinear-form}}{=}&S(D(u_1,\cdots,u_{n-1})(D(v_1,\cdots,v_{n-1})w_1),w_2,\cdots,w_{n-1})\\
&&+\sum_{i=2}^{n-1}S(D(u_1,\cdots,u_{n-1})w_1,w_2,\cdots,w_{i-1},D(v_1,\cdots,v_{n-1})w_i,w_{i+1},\cdots,w_{n-1})\\
&\stackrel{\eqref{unitarity}}{=}&S\Big(\big(D(u_1,\cdots,u_{n-1})\circ D(v_1,\cdots,v_{n-1})-D(v_1,\cdots,v_{n-1})\circ D(u_1,\cdots,u_{n-1})\big)w_1,w_2,\cdots,w_{n-1})\Big)\\
&=&\omega([D(u_1,\cdots,u_{n-1}),D(v_1,\cdots,v_{n-1})]_C,D(w_1,\cdots,w_{n-1})).
\end{eqnarray*}
Therefore, the bilinear form $\omega$ on $\g$ is symmetric, non-degenerate and ad-invariant. The proof is finished. \qed\vspace{3mm}

It is obvious that  $\huaV$ is a faithful representation of the Lie algebra $\g$. Furthermore, we have
\begin{pro}
 $\huaV$ is a faithful generalized orthogonal representation of the Lie algebra $\g$.
\end{pro}
\pf  By the unitarity condition \eqref{unitarity} of a generalized metric $n$-Leibniz algebra, we have
\begin{eqnarray*}
&&S(D(u_1,\cdots,u_{n-1})  w_1,w_2,\cdots,w_{n-1})\\&=&S([u_1,\cdots,u_{n-1},w_1],w_2,\cdots,w_{n-1})\\
                                                   &=&-\sum_{i=2}^{n-1}S(w_1,\cdots,w_{i-1},[u_1,\cdots,u_{n-1},w_i],w_{i+1},\cdots,w_{n-1})\\
                                                   &=&-\sum_{i=2}^{n-1}S(w_1,\cdots,w_{i-1},D(u_1,\cdots,u_{n-1}) w_i,w_{i+1},\cdots,w_{n-1}).
\end{eqnarray*}
Thus, $\huaV$ is faithful generalized orthogonal representation of $\g$. \qed\vspace{3mm}

Summarizing the above discussion, we have

\begin{thm}\label{thm:main1}
  Let $(\huaV,[\cdot,\cdots,\cdot],S)$ be a generalized metric $n$-Leibniz algebra. Then $(\g,\huaV,{\Id})$ is a Lie triple data, i.e. $(\g,\omega)$ is a metric Lie algebra and $({\Id},\huaV,S)$ is its faithful generalized orthogonal representation.
\end{thm}

\begin{ex}{\rm
  Consider the $4$-dimensional simple $3$-Lie algebra on $\mathbb R^4$ with the standard Euclidean structure. With respect to an orthogonal basis $\{e_1,e_2,e_3,e_4\}$, the 3-Lie bracket is given by
 $$[e_1,e_2,e_3]=e_4, \quad [e_2,e_3,e_4]=-e_1,\quad
[e_1,e_3,e_4]=e_2,\quad[e_1,e_2,e_4]=-e_3.$$
It is obvious that $\wedge^2\mathbb R^4$ is 6-dimensional and generated by
$$
e_1\wedge e_2,\quad e_1\wedge e_3,\quad e_1\wedge e_4,\quad e_2\wedge e_3,\quad e_2\wedge e_4,\quad e_3\wedge e_4.
$$
Denote $D(e_i\wedge e_j)$ by $D_{ij}$. We have
$$
D_{12}=\left(\begin{array}{cccc}
  0&0&0&0\\
  0&0&0&0\\
  0&0&0&1\\
  0&0&-1&0
\end{array}\right),~ D_{13}=\left(\begin{array}{cccc}
  0&0&0&0\\
  0&0&0&-1\\
  0&0&0&0\\
  0&1&0&0
\end{array}\right),~ D_{14}=\left(\begin{array}{cccc}
  0&0&0&0\\
  0&0&1&0\\
  0&-1&0&0\\
  0&0&0&0
\end{array}\right),
$$
$$
D_{23}=\left(\begin{array}{cccc}
  0&0&0&1\\
  0&0&0&0\\
  0&0&0&0\\
  -1&0&0&0
\end{array}\right),~ D_{24}=\left(\begin{array}{cccc}
  0&0&-1&0\\
  0&0&0&0\\
  1&0&0&0\\
  0&0&0&0
\end{array}\right),~ D_{34}=\left(\begin{array}{cccc}
  0&1&0&0\\
  -1&0&0&0\\
  0&0&0&0\\
  0&0&0&0
\end{array}\right).
$$
It is obvious that $\{D_{ij}, ~i<j\}$ are basis of $\so(4)$. Therefore, $\g=\Img(D)=\so(4)$.

Next we consider the induced nondegenerate bilinear form $\omega$ on $\g$. The nonzero ones are given by
$$
\omega(D_{12},D_{34})=1,\quad \omega(D_{13},D_{24})=-1, \quad \omega(D_{14},D_{23})=1,
$$
which implies that $\omega$ is not positive definite, but have signature $(3,3)$.}
\end{ex}

\subsection{From the Leibniz algebra to Lie algebra}
Let $(\huaV,[\cdot,\cdots,\cdot],S)$ be a generalized metric $n$-Leibniz algebra. In the middle of the $n$-Leibniz algebra $(\huaV,[\cdot,\cdots,\cdot])$ and the Lie algebra $\g$, we have the Leibniz algebra ($\otimes^{n-1}\huaV, [\cdot,\cdot]_{\HF}$). Moreover, $D$ is a Leibniz algebra epimorphism from $\otimes^{n-1}\huaV$ to $\g$. In this section, we analyze the metric structure on the Leibniz algebra ($\otimes^{n-1}\huaV, [\cdot,\cdot]_{\HF}$). We define a bilinear form $B$ on $\otimes^{n-1}\huaV$ by
\begin{eqnarray}\label{leibniz-bilinear-form}
B(u_1\otimes\cdots \otimes u_{n-1},v_1\otimes\cdots \otimes v_{n-1})=S([u_1,\cdots,u_{n-1},v_1],v_2,\cdots,v_{n-1}).
\end{eqnarray}

\begin{pro}
The bilinear form $B$ on $\otimes^{n-1}\huaV$ defined by \eqref{leibniz-bilinear-form} is symmetric and associative-invariant.
\end{pro}

\pf By the symmetry condition \eqref{symmetry} of a generalized metric $n$-Leibniz algebra, we have
\begin{eqnarray*}
B(u_1\otimes\cdots \otimes u_{n-1},v_1\otimes\cdots \otimes v_{n-1})&=&S([u_1,\cdots,u_{n-1},v_1],v_2,\cdots,v_{n-1})\\
                                                                    &=&S([v_1,v_2,\cdots,v_{n-1},u_1],u_2,\cdots,u_{n-1})\\
                                                                    &=&B(v_1\otimes\cdots \otimes v_{n-1},u_1\otimes\cdots \otimes u_{n-1}).
\end{eqnarray*}
Moreover,  we prove the associative-invariance of the bilinear form $B$:
\begin{eqnarray*}
&&B(u_1\otimes\cdots \otimes u_{n-1},[v_1\otimes\cdots \otimes v_{n-1},w_1\otimes\cdots \otimes w_{n-1}]_{\HF})\\
&\stackrel{\eqref{leibniz-bracket}}{=}&B(u_1\otimes\cdots \otimes u_{n-1},\sum_{i=1}^{n-1}w_1\otimes\cdots w_{i-1}\otimes [v_1,\cdots,v_{n-1},w_i]\otimes w_{i+1}\cdots \otimes w_{n-1}))\\
&\stackrel{\eqref{leibniz-bilinear-form}}{=}&S([u_1,\cdots,u_{n-1},[v_1,\cdots,v_{n-1},w_1]],w_2,\cdots,w_{n-1})\\
&&+\sum_{i=2}^{n-1}S([u_1,\cdots,u_{n-1},w_1],w_2,\cdots,w_{i-1},[v_1,\cdots,v_{n-1},w_i],w_{i+1},\cdots,w_{n-1})\\
&\stackrel{\eqref{unitarity}}{=}&S([u_1,\cdots,u_{n-1},[v_1,\cdots,v_{n-1},w_1]],w_2,\cdots,w_{n-1})\\
&&-S([v_1,\cdots,v_{n-1},[u_1,\cdots,u_{n-1},w_1]],w_2,\cdots,w_{n-1})\\
&\stackrel{\eqref{fundamental-identity}}{=}&B([u_1\otimes\cdots \otimes u_{n-1},v_1\otimes\cdots \otimes v_{n-1}]_{\HF},w_1\otimes\cdots \otimes w_{n-1}).
\end{eqnarray*}
Therefore, the bilinear form $B$ on $\otimes^{n-1}\huaV$ is symmetric and associative-invariant. The proof is finished. \qed

\begin{pro}
The bilinear form $B$ on $\otimes^{n-1}\huaV$ is non-degenerate if and only if $\ker D=0.$
\end{pro}

\pf Let $V=\sum_iv_{i,1}\otimes\cdots v_{i,n-1}\in\otimes^{n-1}\huaV$ be such that $B(V,w_1\otimes\cdots\otimes w_{n-1})=0$ for all $w_1,w_2,\cdots,w_{n-1}\in\huaV$. Therefore we have
\begin{eqnarray*}
S([V,w_1],w_2,\cdots,w_{n-1})=0.
\end{eqnarray*}
Since $S$ is non-degenerate, we have $[V,w_1]=0$ for all $w_1\in\huaV$, hence $V\in\ker D$.
The proof is finished. \qed

\begin{rmk}
Since the Leibniz algebra $(\otimes^{n-1}\huaV, [\cdot,\cdot]_{\HF})$ is not an anticommutative algebra in general. Thus, a symmetric  associative-invariant bilinear form $B$ is not $\ad$-invariant.
\end{rmk}

\begin{pro}
The Leibniz algebra morphism $D:\otimes^{n-1}\huaV\longrightarrow\g$ preserves the metric.
\end{pro}

\pf For all $u_1,\cdots,u_{n-1},v_1,\cdots,v_{n-1}\in\huaV$, we have
\begin{eqnarray*}
\omega(D(u_1,\cdots,u_{n-1}),D(v_1,\cdots,v_{n-1}))&=&S([u_1,\cdots,u_{n-1},v_1],v_2\cdots, v_{n-1})\\
                                             &=&B(u_1\otimes\cdots \otimes u_{n-1},v_1\otimes\cdots \otimes v_{n-1}).
\end{eqnarray*}
Thus, $D$ preserves the metric. \qed

\section{Construction of a generalized metric $n$-Leibniz algebra from a Lie triple data}
Let $(\g,[\cdot,\cdot],\omega)$ be a metric  Lie algebra and $(\rho,\huaV,S)$ a faithful generalized orthogonal representation of $\g$ as defined in Definition \ref{orthogonal-representation}. We start by defining an $(n-1)$-linear map $D:\huaV\times\cdots\times\huaV\lon\g$, by transposing the $\g$-action. That is, for given $v_1,\cdots,v_{n-1}\in\huaV$, define $D(v_1,\cdots,v_{n-1})\in\g$  by
\begin{eqnarray}\label{Faulkner-construction}
\omega(x,D(v_1,\cdots,v_{n-1}))=S(\rho(x)v_1,v_2,\cdots,v_{n-1}),\,\,\,\,\forall x\in\g.
\end{eqnarray}

\begin{pro}\label{skew-symmetry}
With the above notations, for all $v_1,v_2,\cdots,v_{n-1}\in\huaV$, we have $$\sum_{i=1}^{n-1}D(v_i,v_1,\cdots,v_{i-1},\hat{v_i},v_{i+1},\cdots,v_{n-1})=0.$$
\end{pro}

\pf Since $(\rho,\huaV,S)$ is a generalized orthogonal representation of $\g$, we have
\begin{eqnarray*}
\omega(x,D(v_1,\cdots,v_{n-1}))&=&S(\rho(x)v_1,v_2,\cdots,v_{n-1})\\
                              &\stackrel{\eqref{orthogonal-representation-invariant}}{=}&-\sum_{i=2}^{n-1}S(v_1,\cdots,v_{i-1},\rho(x)v_i,v_{i+1},\cdots,v_{n-1})\\
                              &=&-\sum_{i=2}^{n-1}S(\rho(x)v_i,v_1,\cdots,v_{i-1},\hat{v_i},v_{i+1},\cdots,v_{n-1})\\
                              &=&-\sum_{i=2}^{n-1}\omega(x,D(v_i,v_1,\cdots,v_{i-1},\hat{v_i},v_{i+1},\cdots,v_{n-1})).
\end{eqnarray*}
By the nondegeneracy of $\omega$, we have $$D(v_1,v_2,\cdots,v_{n-1})=-\sum_{i=2}^{n-1}D(v_i,v_1,\cdots,v_{i-1},\hat{v_i},v_{i+1},\cdots,v_{n-1}).$$
Thus, the proof is finished. \qed

\begin{pro}
The $(n-1)$-linear map $D:\huaV\times\cdots\times\huaV\lon\g$ is  surjective.
\end{pro}

\pf We denote by $(\Img D)^{\bot}$ the orthogonal compliment space of $\Img D$, i.e.
$$
(\Img D)^{\bot}:=\{x\in\g|\omega(x,y)=0,\,\, \forall y\in\Img D\}.
$$
Let $x\in (\Img D)^{\bot}$.
Then for all $v_1,\cdots,v_{n-1}\in\huaV$, we have $$\omega(x,D(v_1,\cdots,v_{n-1}))=0.$$ Therefore, by \eqref{Faulkner-construction} we obtain $S(\rho(x) v_1,v_2,\cdots,v_{n-1})=0$. The nondegeneracy of $S$ implies that $\rho(x) v_1=0$ for all $v_1\in\huaV$, which in turn implies that $x=0$ since the representation of $\g$ on $\huaV$ is faithful. Therefore, $(\Img D)^{\bot}=0$ and $D$ is surjective.   \qed

\vspace{2mm}
We define an $n$-linear map $[\cdot,\cdots,\cdot]:\huaV\times\cdots\times\huaV\lon\huaV$ by
\begin{eqnarray}\label{n-bracket}
[v_1,\cdots,v_{n-1},v_n]=\rho(D(v_1,\cdots,v_{n-1}))  v_n.
\end{eqnarray}
By  Proposition \ref{skew-symmetry}, it is straightforward to obtain
\begin{lem}
For all $v_1,\cdots,v_{n}\in\huaV$, there holds
\begin{eqnarray*}
\sum_{i=1}^{n-1}[v_i,v_1,\cdots,v_{i-1},\hat{v_i},v_{i+1},\cdots,v_{n-1},v_n]=0.
\end{eqnarray*}
\end{lem}

\begin{rmk}
For $n=3$, we obtain that the $3$-bracket is skew-symmetric in the first two entries.
\end{rmk}

\vspace{2mm}
The following theorem says that the converse of Theorem \ref{thm:main1} also holds. Thus, there is a one-to-one correspondence between generalized metric $n$-Leibniz algebras and faithful generalized orthogonal representations of metric Lie algebras.

\begin{thm}\label{thm:main2}
Let  $(\rho,\huaV,S)$ be a faithful generalized orthogonal representation of a metric  Lie algebra $(\g,[\cdot,\cdot],\omega)$. Then $(\huaV,[\cdot,\cdots,\cdot],S)$ is a generalized metric $n$-Leibniz algebra, where the $n$-bracket $[\cdot,\cdots,\cdot]$ is defined by \eqref{n-bracket}.
\end{thm}

\pf For all $x,y\in\g$ and $u_1,\cdots,u_{n-1}\in\huaV$, we have
\begin{eqnarray*}
\omega([D(u_1,\cdots,u_{n-1}),x],y)&\stackrel{\eqref{ad-invariant}}{=}&\omega(D(u_1,\cdots,u_{n-1}),[x,y])\\
                                           &=&\omega([x,y],D(u_1,\cdots,u_{n-1}))\\
                                           &\stackrel{\eqref{Faulkner-construction}}{=}&S(\rho([x,y]) u_1,u_2,\cdots,u_{n-1})\\
                                           &=&S(\rho(x)\rho(y) u_1,u_2,\cdots,u_{n-1})-S(\rho(y)\rho(x) u_1,u_2,\cdots,u_{n-1})\\
                                          &\stackrel{\eqref{orthogonal-representation-invariant}}{=}&-\sum_{i=2}^{n-1}S(\rho(y) u_1,u_2,\cdots,u_{i-1},\rho(x) u_i,u_{i+1},\cdots,u_{n-1})\\
                                          &&-S(\rho(y)\rho(x) u_1,u_2,\cdots,u_{n-1})\\
                                          &\stackrel{\eqref{Faulkner-construction}}{=}&-\sum_{i=2}^{n-1}\omega(y,D(u_1,u_2,\cdots,u_{i-1},\rho(x) u_i,u_{i+1},\cdots,u_{n-1}))\\
                                          &&-\omega(y,D(\rho(x) u_1,u_2,\cdots,u_{n-1})).
\end{eqnarray*}
By the nondegeneracy of the bilinear form $\omega$ on $\g$, we have
\begin{eqnarray*}
[x,D(u_1,\cdots,u_{n-1})]=\sum_{i=1}^{n-1}D(u_1,\cdots,u_{i-1},\rho(x) u_i,u_{i+1},\cdots,u_{n-1}).
\end{eqnarray*}
By substituting $x=D(v_1,\cdots,v_{n-1})$ and applying both sides of the above equation to $u_n$, we have
\begin{eqnarray*}
&&[v_1,\cdots,v_{n-1},[u_1,\cdots,u_{n-1},u_n]]-[u_1,\cdots,u_{n-1},[v_1,\cdots,v_{n-1},u_n]]\\
&=&\sum_{i=1}^{n-1}[u_1,\cdots,u_{i-1},[v_1,\cdots,v_{n-1},u_i],u_{i+1},\cdots,u_{n-1},u_n].
\end{eqnarray*}
Thus, $(\huaV,[\cdot,\cdots,\cdot])$ is an $n$-Leibniz algebra.

By \eqref{orthogonal-representation-invariant}, we have
\begin{eqnarray*}
&&\sum_{i=1}^{n-1}S(v_1,\cdots,v_{i-1},[u_1,\cdots,u_{n-1}, v_i],v_{i+1},\cdots v_{n-1})\\&=& \sum_{i=1}^{n-1}S(v_1,\cdots,v_{i-1},\rho(D(u_1,\cdots,u_{n-1})) v_i,v_{i+1},\cdots v_{n-1})\\
                                               &=&0.
\end{eqnarray*}
Thus, the unitarity condition in Definition \ref{n-leibniz-algebra} holds.

Since the bilinear form $\omega$ on $\g$ is symmetric, we have
\begin{eqnarray*}
S([u_1,\cdots,u_{n-1},v_1],v_2,\cdots, v_{n-1})&=&S(\rho(D(u_1,\cdots,u_{n-1})) v_1,v_2,\cdots, v_{n-1})\\
                                               &\stackrel{\eqref{Faulkner-construction}}{=}&\omega(D(u_1,\cdots,u_{n-1}),D(v_1,v_2,\cdots, v_{n-1}))\\
                                               &=&\omega(D(v_1,v_2,\cdots, v_{n-1}),D(u_1,\cdots,u_{n-1}))\\
                                               &\stackrel{\eqref{Faulkner-construction}}{=}&S(\rho(D(v_1,v_2,\cdots, v_{n-1})) u_1,u_2,\cdots,u_{n-1})\\
                                               &=&S([v_1,v_2,\cdots, v_{n-1},u_1],u_2,\cdots,u_{n-1}),
\end{eqnarray*}
which implies that the symmetry condition in Definition \ref{n-leibniz-algebra} holds.

Thus, $(\huaV,[\cdot,\cdots,\cdot],S)$ is a generalized   metric $n$-Leibniz algebra. The proof is finished. \qed

\emptycomment{
\subsection{Isomorphic of two groupoids  }
\begin{defi}
An isomorphism of generalized metric $n$-Leibniz algebras $\varphi:(\huaV,[\cdot,\cdots,\cdot]_{\huaV},S_{\huaV})\lon(\huaW,[\cdot,\cdots,\cdot]_{\huaW},S_{\huaW})$ is a linear isomorphism $\varphi:\huaV\lon\huaW$ such that
\begin{eqnarray}
\varphi[v_1,\cdots,v_{n}]_{\huaV}&=&[\varphi(v_1),\cdots,\varphi(v_{n})]_{\huaW}\\
S_{\huaV}(v_1,\cdots,v_{n-1})&=&S_{\huaW}(\varphi(v_1),\cdots,\varphi(v_{n-1})).
\end{eqnarray}
\end{defi}

We denote by {\bf generalized metric $n$-Leibniz} the category of generalized metric $n$-Leibniz algebras and their isomorphisms. Since every morphism is invertible, we obtain that {\bf generalized metric $n$-Leibniz} is a groupoid.

\begin{defi}
An isomorphism of the pairs of a metric Lie algebra and its a faithful generalized orthogonal representation $f:(\g_{\huaV}\subset\gl(\huaV))\lon(\g_{\huaW}\subset\gl(\huaW))$ is a linear isomorphism $f:\huaV\lon\huaW$ such that $\Ad_{f}:\g_{\huaV}\lon\g_{\huaW}$ is an isometry of metric Lie algebras and
\begin{eqnarray}
S_{\huaV}(v_1,\cdots,v_{n-1})&=&S_{\huaW}(f(v_1),\cdots,f(v_{n-1})).
\end{eqnarray}
\end{defi}

We denote by {\bf pairs } the category of the pairs $(\g_{\huaV}\subset\gl(\huaV))$ and their isomorphisms. Since every morphism is invertible, we obtain that {\bf pairs} is a groupoid.

\begin{thm}
The two groupoids are isomorphic.
\end{thm}
\pf Let $(\huaV,[\cdot,\cdots,\cdot]_{\huaV},S_{\huaV})$ be a generalized metric $n$-Leibniz algebra
}

\section{Generalized  orthogonal derivations}

In this section, we introduce the notion of a generalized orthogonal derivation on a Lie triple data and show that there is a one-to-one correspondence between generalized orthogonal derivations on generalized metric $n$-Leibniz algebras and Lie triple datas.

\begin{defi} A generalized orthogonal derivation on a Lie triple data $(\g,\huaV,\rho)$ is a pair $(\dM_\g,\dM_\huaV)$, where $\dM_{\g}$ is an orthogonal derivation on the metric Lie algebra $(\g,[\cdot,\cdot],\omega)$ and $\dM_{\huaV}\in\gl(\huaV)$ is a linear map satisfying the following conditions:
\begin{eqnarray}
\label{orthogonal-derivation-1}&&\dM_{\huaV}\circ\rho(x)=\rho(\dM_{\g}(x))+\rho(x)\circ\dM_{\huaV},\\
\label{orthogonal-derivation-2}&&\sum_{i=1}^{n-1}S(w_1,\cdots,\dM_{\huaV}w_i,\cdots,w_{n-1})=0,
\end{eqnarray}
for all $x\in\g$ and $w_1,w_2,\cdots,w_{n-1}\in \huaV$.
\end{defi}

Let $(\huaV,[\cdot,\cdots,\cdot],S)$ be a generalized metric $n$-Leibniz algebra with a generalized  orthogonal derivation $\dM_{\huaV}$.  Let $(\g,[\cdot,\cdot]_C,\omega)$ be the corresponding metric Lie algebra given in Proposition \ref{pro:metLie}.
Define $\dM_{\g}:\g\longrightarrow \g$   by
\begin{eqnarray}
\dM_{\g}\big(D(w_1,\cdots,w_{n-1})\big)=\sum_{i=1}^{n-1}D(w_1,\cdots,\dM_{\huaV}w_i,\cdots,w_{n-1}).
\end{eqnarray}
Or equivalently,
\begin{eqnarray*}
\dM_{\g}\big(D(w_1,\cdots,w_{n-1})\big)=[\dM_{\huaV},D(w_1,\cdots,w_{n-1})]_C.
\end{eqnarray*}

\begin{pro}
Let $\dM_{\huaV}$ be a generalized orthogonal derivation on a generalized metric $n$-Leibniz algebra $(\huaV,[\cdot,\cdots,\cdot],S)$.  Then $(\dM_{\g},\dM_\huaV)$ is a generalized  orthogonal derivation on the Lie triple data $(\g,\huaV,\Id)$ given by Theorem \ref{thm:main1}.
\end{pro}
\pf For all  $u_1,\cdots,u_{n-1},v_1,\cdots,v_{n-1}\in \huaV$, we have
\begin{eqnarray*}
 && \dM_\g[D(u_1,\cdots,u_{n-1}),D(v_1,\cdots,v_{n-1})]_C\\&=&[\dM_\huaV,[D(u_1,\cdots,u_{n-1}),D(v_1,\cdots,v_{n-1})]_C]_C\\
  &=&[[\dM_\huaV,D(u_1,\cdots,u_{n-1})]_C,D(v_1,\cdots,v_{n-1})]_C+[D(u_1,\cdots,u_{n-1}),[\dM_\huaV,D(v_1,\cdots,v_{n-1})]_C]_C\\
  &=&[\dM_\g\big(D(u_1,\cdots,u_{n-1})\big),D(v_1,\cdots,v_{n-1})]_C+[D(u_1,\cdots,u_{n-1}),\dM_\g\big(D(v_1,\cdots,v_{n-1})\big)]_C,
\end{eqnarray*}
which implies that $\dM_\g$ is a derivation of the Lie algebra $(\g,[\cdot,\cdot]_C)$.

Since  $\dM_{\huaV}$ is generalized orthogonal, for all $D(u_1,\cdots,u_{n-1}),D(v_1,\cdots,v_{n-1})\in\g$, we have
\begin{eqnarray*}
&&\omega(\dM_{\g}D(u_1,\cdots,u_{n-1}),D(v_1,\cdots,v_{n-1}))+\omega(D(u_1,\cdots,u_{n-1}),\dM_{\g}D(v_1,\cdots,v_{n-1}))\\
&=&\sum_{i=1}^{n-1}S([u_1,\cdots,\dM_{\huaV}u_i,\cdots,u_{n-1},v_1],v_2,\cdots,v_{n-1})+S([u_1,\cdots,u_{n-1},\dM_{\huaV}v_1],v_2,\cdots,v_{n-1})\\
&&+\sum_{i=2}^{n-1}S([u_1,\cdots,u_{n-1},v_1],v_2,\cdots,\dM_{\huaV}v_i,\cdots,v_{n-1})\\
&=&S(\dM_{\huaV}[u_1,\cdots,u_{n-1},v_1],v_2,\cdots,v_{n-1})+\sum_{i=2}^{n-1}S([u_1,\cdots,u_{n-1},v_1],v_2,\cdots,\dM_{\huaV}v_i,\cdots,v_{n-1})\\
&=&0.
\end{eqnarray*}
Thus, $\dM_{\g}$ is an orthogonal derivation on the metric Lie algebra $(\g,[\cdot,\cdot]_C,\omega)$.

Moreover, for all $D(u_1,\cdots,u_{n-1})\in\g$, we have
\begin{eqnarray*}
\dM_{\g}(D(u_1,\cdots,u_{n-1}))+D(u_1,\cdots,u_{n-1})\circ\dM_{\huaV}
&=&[\dM_{\huaV},D(u_1,\cdots,u_{n-1})]_C+D(u_1,\cdots,u_{n-1})\circ\dM_{\huaV}\\
&=&\dM_{\huaV}\circ D(u_1,\cdots,u_{n-1}).
\end{eqnarray*}
Thus, equality \eqref{orthogonal-derivation-1} holds. Furthermore,   \eqref{orthogonal-derivation-2} holds automatically. The proof is finished. \qed\vspace{3mm}

The converse of the above result also holds.

\begin{pro}
Let $(\dM_\g,\dM_{\huaV})$ be a generalized orthogonal derivation on a Lie triple data $(\g,\huaV,\rho)$.  Then $\dM_{\huaV}$ is a generalized orthogonal derivation on the corresponding generalized metric $n$-Leibniz algebra $(\huaV,[\cdot,\cdots,\cdot],S)$ given in Theorem \ref{thm:main2}.
\end{pro}
\pf  We only need to prove
that $\dM_{\huaV}$ is a derivation on the $n$-Leibniz algebra $(\huaV,[\cdot,\cdots,\cdot])$.
For all $v_1,\cdots,v_{n-1}\in\huaV$ and $x\in\g$, we have
\begin{eqnarray*}
&&\omega(\dM_{\g}D(v_1,\cdots,v_{n-1})-\sum_{i=1}^{n-1}D(v_1,\cdots,\dM_{\huaV}v_i,\cdots,v_{n-1}),x)\\
&\stackrel{\eqref{o-derivation}}{=}&-\omega(D(v_1,\cdots,v_{n-1}),\dM_{\g}x)-\sum_{i=2}^{n-1}S(\rho(x)v_1,v_2,\cdots,\dM_{\huaV}v_i,\cdots,v_{n-1})-
S(\rho(x)(\dM_{\huaV}v_1),v_2,\cdots,v_{n-1})\\
&\stackrel{\eqref{orthogonal-derivation-2}}{=}&-S(\rho(\dM_{\g}x)v_1,v_2,\cdots,v_{n-1})+S(\dM_{\huaV}(\rho(x)v_1),v_2,\cdots,v_i,\cdots,v_{n-1})\\&&-
S(\rho(x)(\dM_{\huaV}v_1),v_2,\cdots,v_{n-1})\\
&\stackrel{\eqref{orthogonal-derivation-1}}{=}&0.
\end{eqnarray*}
Thus, we have
\begin{eqnarray}\label{derivation-identity}
\dM_{\g}D(v_1,\cdots,v_{n-1})=\sum_{i=1}^{n-1}D(v_1,\cdots,\dM_{\huaV}v_i,\cdots,v_{n-1}).
\end{eqnarray}
For all $v_1,\cdots,v_n,u_1,\cdots,u_{n-2}\in\huaV$, we have
\begin{eqnarray*}
&&S(\dM_{\huaV}[v_1,\cdots,v_n]-\sum_{i=1}^{n}[v_1,\cdots,\dM_{\huaV}v_i,\cdots,v_n],u_1,\cdots,u_{n-2})\\
&=&S(\dM_{\huaV}(\rho(D(v_1,\cdots,v_{n-1}))v_n),u_1,\cdots,u_{n-2})\\
&&-\sum_{i=1}^{n-1}S(\rho(D(v_1,\cdots,\dM_{\huaV}v_i,\cdots,v_{n-1}))v_n,u_1,\cdots,u_{n-2})\\
&&-S(\rho(D(v_1,\cdots,v_{n-1}))(\dM_{\huaV}v_n),u_1,\cdots,u_{n-2})\\
&\stackrel{\eqref{orthogonal-derivation-1}}{=}&S(\rho(\dM_{\g}D(v_1,\cdots,v_{n-1}))v_n,u_1,\cdots,u_{n-2})\\
&&-\sum_{i=1}^{n-1}S(\rho(D(v_1,\cdots,\dM_{\huaV}v_i,\cdots,v_{n-1}))v_n,u_1,\cdots,u_{n-2})\\
&\stackrel{\eqref{derivation-identity}}{=}&0.
\end{eqnarray*}
Therefore,  $\dM_{\huaV}$ is a derivation of the $n$-Leibniz algebra $(\huaV,[\cdot,\cdots,\cdot])$. \qed

\section{Generalized  orthogonal automorphisms}

In this section, we introduce the notion of a generalized orthogonal automorphism on a Lie triple data and show that there is a one-to-one correspondence between generalized orthogonal automorphisms on generalized metric $n$-Leibniz algebras and Lie triple datas.

\begin{defi} A generalized orthogonal automorphism on a Lie triple data $(\g,\huaV,\rho)$ is a pair $(\Phi_\g,\Phi_\huaV)$, where $\Phi_\g$ is an orthogonal automorphism on the metric Lie algebra $(\g,[\cdot,\cdot],\omega)$ and $\Phi_\huaV\in\gl(\huaV)$ is an invertible linear map satisfying the following conditions:
\begin{eqnarray}
\label{orthogonal-automorphism-1}&&\Phi_\huaV(\rho(x)w)=\rho(\Phi_\g(x))(\Phi_\huaV w),\\
\label{orthogonal-automorphism-2}&&S(\Phi_\huaV w_1,\cdots,\Phi_\huaV w_{n-1})=S( w_1,\cdots, w_{n-1}),
\end{eqnarray}
for all $x\in\g$ and $w,w_1,w_2,\cdots,w_{n-1}\in \huaV$.
\end{defi}

Let $(\huaV,[\cdot,\cdots,\cdot],S)$ be a generalized metric $n$-Leibniz algebra with a generalized  orthogonal automorphism $\Phi_{\huaV}$.  Let $(\g,[\cdot,\cdot]_C,\omega)$ be the corresponding metric Lie algebra given in Proposition \ref{pro:metLie}.
Define $\Phi_{\g}:\g\longrightarrow \g$   by
\begin{eqnarray}
\Phi_{\g}\big(D(w_1,\cdots,w_{n-1})\big)=D(\Phi_{\huaV}w_1,\cdots,\Phi_{\huaV}w_{n-1}).
\end{eqnarray}
Or equivalently,
\begin{eqnarray*}
\Phi_{\g}\big(D(w_1,\cdots,w_{n-1})\big)=\Phi_{\huaV}\circ D(w_1,\cdots,w_{n-1})\circ\Phi_{\huaV}^{-1}.
\end{eqnarray*}

\begin{pro}
Let $\Phi_{\huaV}$ be an generalized orthogonal automorphism on a generalized metric $n$-Leibniz algebra $(\huaV,[\cdot,\cdots,\cdot],S)$.  Then $(\Phi_\g,\Phi_\huaV)$ is a generalized  orthogonal automorphism on the Lie triple data $(\g,\huaV,\Id)$ given by Theorem \ref{thm:main1}.
\end{pro}
\pf For all  $u_1,\cdots,u_{n-1},v_1,\cdots,v_{n-1}\in \huaV$, we have
\begin{eqnarray*}
&&\Phi_{\g}[D(u_1,\cdots,u_{n-1}),D(v_1,\cdots,v_{n-1})]_C\\
&=&\Phi_{\huaV}\circ D(u_1,\cdots,u_{n-1})\circ D(v_1,\cdots,v_{n-1})\circ\Phi_{\huaV}^{-1}\\&&-\Phi_{\huaV}\circ D(v_1,\cdots,v_{n-1})\circ D(u_1,\cdots,u_{n-1})\circ\Phi_{\huaV}^{-1}\\
&=&[\Phi_{\g}D(u_1,\cdots,u_{n-1}),\Phi_{\g}D(v_1,\cdots,v_{n-1})]_C.
\end{eqnarray*}
Thus, $\Phi_{\g}$ is an automorphism of the Lie algebra $(\g,[\cdot,\cdot]_C)$. Since  $\Phi_{\huaV}$ is generalized orthogonal, for all $D(u_1,\cdots,u_{n-1}),D(v_1,\cdots,v_{n-1})\in\g$, we have
\begin{eqnarray*}
&&\omega(\Phi_{\g}D(u_1,\cdots,u_{n-1}),\Phi_{\g}D(v_1,\cdots,v_{n-1}))\\
&=&S([\Phi_{\huaV}u_1,\cdots,\Phi_{\huaV}u_{n-1},\Phi_{\huaV}v_1],\Phi_{\huaV}v_2,\cdots,\Phi_{\huaV}v_{n-1})\\
&=&S(\Phi_{\huaV}[u_1,\cdots,u_{n-1},v_1],\Phi_{\huaV}v_2,\cdots,\Phi_{\huaV}v_{n-1})\\
&=&S([u_1,\cdots,u_{n-1},v_1],v_2,\cdots,v_{n-1})\\
&=&\omega(D(u_1,\cdots,u_{n-1}),D(v_1,\cdots,v_{n-1})).
\end{eqnarray*}
Thus, $\Phi_{\g}$ is an orthogonal automorphism on the metric Lie algebra $(\g,[\cdot,\cdot]_C,\omega)$.

Moreover, for all $D(u_1,\cdots,u_{n-1})\in\g$ and $w\in \huaV$, we have
\begin{eqnarray*}
\Phi_\huaV(D(u_1,\cdots,u_{n-1})w)&=&[\Phi_\huaV u_1,\cdots,\Phi_\huaV u_{n-1},\Phi_\huaV w]\\
                                  &=&(\Phi_{\g}D(u_1,\cdots,u_{n-1}))(\Phi_\huaV w).
\end{eqnarray*}
Thus, equality \eqref{orthogonal-automorphism-1} holds. Furthermore,   \eqref{orthogonal-automorphism-2} holds automatically. The proof is finished. \qed\vspace{3mm}

The converse of the above result also holds.
\begin{pro}
Let $(\Phi_\g,\Phi_{\huaV})$ be a generalized orthogonal automorphism on a Lie triple data $(\g,\huaV,\rho)$. Then $\Phi_{\huaV}$ is a generalized orthogonal automorphism on the corresponding generalized metric $n$-Leibniz algebra $(\huaV,[\cdot,\cdots,\cdot],S)$ given in Theorem \ref{thm:main2}.
\end{pro}
\pf  We only need to prove
that $\Phi_{\huaV}$ is an automorphism on the $n$-Leibniz algebra $(\huaV,[\cdot,\cdots,\cdot])$.
For all $v_1,\cdots,v_{n-1}\in\huaV$ and $x\in\g$, we have
\begin{eqnarray*}
&&\omega(\Phi_{\g}D(v_1,\cdots,v_{n-1})-D(\Phi_{\huaV}v_1,\cdots,\Phi_{\huaV}v_{n-1}),x)\\
&\stackrel{\eqref{o-automorphism}}{=}&\omega(D(v_1,\cdots,v_{n-1}),\Phi_{\g}^{-1}x)-\omega(D(\Phi_{\huaV}v_1,\cdots,\Phi_{\huaV}v_{n-1}),x)\\
&=&S(\rho(\Phi_{\g}^{-1}x)v_1,v_2,\cdots,v_{n-1})-S(\rho(x)(\Phi_{\huaV}v_1),\Phi_{\huaV}v_2,\cdots,\Phi_{\huaV}v_{n-1})\\
&\stackrel{\eqref{orthogonal-automorphism-2}}{=}&S(\Phi_{\huaV}(\rho(\Phi_{\g}^{-1}x)v_1),\Phi_{\huaV}v_2,\cdots,\Phi_{\huaV}v_{n-1})-S(\rho(x)(\Phi_{\huaV}v_1),\Phi_{\huaV}v_2,\cdots,\Phi_{\huaV}v_{n-1})\\
&\stackrel{\eqref{orthogonal-automorphism-1}}{=}&S(\rho(x)(\Phi_{\huaV}v_1),\Phi_{\huaV}v_2,\cdots,\Phi_{\huaV}v_{n-1})-S(\rho(x)(\Phi_{\huaV}v_1),\Phi_{\huaV}v_2,\cdots,\Phi_{\huaV}v_{n-1})\\
&=&0.
\end{eqnarray*}
Thus, we have
\begin{eqnarray}\label{automorphism-identity}
\Phi_{\g}D(v_1,\cdots,v_{n-1})=D(\Phi_{\huaV}v_1,\cdots,\Phi_{\huaV}v_{n-1}).
\end{eqnarray}
For all $v_1,\cdots,v_n,u_1,\cdots,u_{n-2}\in\huaV$, we have
\begin{eqnarray*}
&&S(\Phi_{\huaV}[v_1,\cdots,v_n]-[\Phi_{\huaV}v_1,\cdots,\Phi_{\huaV}v_i,\cdots,\Phi_{\huaV}v_n],u_1,\cdots,u_{n-2})\\
&=&S(\Phi_{\huaV}(\rho(D(v_1,\cdots,v_{n-1}))v_n),u_1,\cdots,u_{n-2})\\&&-S(\rho(D(\Phi_{\huaV}v_1,\cdots,\Phi_{\huaV}v_i,\cdots,\Phi_{\huaV}v_{n-1}))\Phi_{\huaV}v_n,u_1,\cdots,u_{n-2})\\
&\stackrel{\eqref{orthogonal-automorphism-1}}{=}&S(\rho(\Phi_{\g}D(v_1,\cdots,v_{n-1}))(\Phi_{\huaV}v_n),u_1,\cdots,u_{n-2})\\&&-S(\rho(D(\Phi_{\huaV}v_1,\cdots,\Phi_{\huaV}v_i,\cdots,\Phi_{\huaV}v_{n-1}))(\Phi_{\huaV}v_n),u_1,\cdots,u_{n-2})\\
&\stackrel{\eqref{automorphism-identity}}{=}&0.
\end{eqnarray*}
Therefore,  $\Phi_{\huaV}$ is an automorphism of the $n$-Leibniz algebra $(\huaV,[\cdot,\cdots,\cdot])$. \qed

\emptycomment{
\begin{rmk}
When $n=2$, a generalized metric $2$-Leibniz algebra is a Leibniz algebra $(\huaV,[\cdot,\cdot])$ equipped with a linear function $S\in\huaV^*$ such that for all $u,v\in\huaV$, we have $S([u,v])=0$.
\end{rmk}
}

\end{document}